\documentclass[oneside]{amsart}

\usepackage[letterpaper,body={12.6cm,22.5cm}, mag=1000]{geometry}
\usepackage{amssymb}
\usepackage{amsthm}
\usepackage{amscd}

\numberwithin{equation}{section}
\theoremstyle{plain}
\newtheorem{cor}[equation]{Corollary}

\newtheorem{prop}[equation]{Proposition}

\newtheorem{thm}[equation]{Theorem}
\newtheorem{problem}[equation]{Problem}

\theoremstyle{definition}

\newtheorem{remark}[equation]{Remark}

\newcommand{\dlabel}[1]{\ifmmode \text{\ttfamily \upshape [#1] } \else
{\ttfamily \upshape [#1] }\fi \label{#1}}

\newcommand{\C}{\operatorname{C} }
\newcommand{\Z}{\operatorname{Z} }

\newcommand{\gen}[1]{\left < #1 \right >}
\newcommand{\Aut}{\operatorname{Aut} }

\newcommand{\Hom}{\operatorname{Hom} }
\newcommand{\Inn}{\operatorname{Inn} }

\newcommand{\Out}{\operatorname{Out} }

\newcommand{\Autcent}{\operatorname{Autcent} }

\begin{document}

\title{Class preserving automorphisms of finite $p$-groups: A survey}

\author{Manoj K.~Yadav}

\address{School of Mathematics, Harish-Chandra Research Institute \\
Chhatnag Road, Jhunsi, Allahabad - 211 019, INDIA}

\email{myadav@hri.res.in}
\thanks{2000 Mathematics Subject Classification. 20D45, 20D15}
\date{\today}

\maketitle

\begin{abstract}
In this short survey article, we try to list maximum number of known results on class preserving automorphisms of finite $p$-groups.
\end{abstract}

\section{Introduction}

Let us start with a finite group $G$. For $x \in G$, $x^G$ denotes the conjugacy class of 
$x$ in $G$. By $\Aut(G)$ we denote the group of all automorphisms of $G$.
An automorphism $\alpha$ of $G$ is called \emph{class preserving} if 
$\alpha(x) \in x^G$ for all $x \in G$. The set of all class preserving
automorphisms of $G$,  denoted by $\Aut_{c}(G)$, is a normal subgroup of 
$\Aut(G)$. Notice that $\Inn(G)$, the group of all inner automorphisms of $G$,
is a normal subgroup of $\Aut_c(G)$. The group of all class preserving outer automorphisms of $G$, i.e., $\Aut_c(G)/\Inn(G)$ is denoted by $\Out_c(G)$.

The story begins with the following question of W. Burnside \cite[Note B]{wB11}: 
Does there exist any  finite group $G$ such that $G$ has a non-inner 
class preserving automorphism? 
In 1913,  Burnside \cite{wB13} himself gave an affirmative answer to this 
question.
He constructed a group $G$ of order $p^6$ isomorphic to the group $W$ 
consisting  of all $3 \times 3$ matrices
\[  M = \begin{pmatrix} 
          1 & 0 & 0 \\
          x & 1 & 0 \\
          z & y & 1
        \end{pmatrix}  \]
with $x,y,z$ in the field $\mathbb F_{p^2}$ of $p^2$ elements, where $p$ is an
odd prime.
For this group $G$,  $\Out_{c}(G) \not= 1$.
He also proved that $\Aut_{c}(G)$ is an elementary abelian $p$-group of order
$p^8$.

In 1947, G. E. Wall \cite{gW47} constructed examples of arbitrary finite groups $G$ such that $\Out_c(G) \not= 1$. Interestingly his examples contain $2$-group having class preserving outer automorphisms. The smallest of these groups is a group of order $32$. The members of the class of groups, constructed by Wall, appear as a semidirect product of a cyclic group by an abelian group. 

At that stage we had examples of arbitrary finite groups having class preserving outer automorphims. But, unfortunately, this problem was not taken up (upto my knowledge) for next 20 years. However in this span of time many results were proved on automorphisms of finite groups, and more precisely on automorphisms of finite $p$-groups. But all of these results were centered around finding a non-inner $p$-automorphism (an automorphism having order a power of a prime $p$) of finite $p$-groups with order at least $p^2$. This was finally achieved by W. Gaschutz in 1966 \cite{wG66} by using cohomology of finite groups. 

In 1968, C. H. Sah \cite{cS68} first time studied some properties of $\Aut_c(G)$ for an arbitrary group $G$. He used the concepts of $c$-closed and strongly $c$-closed subgroups (cf. \cite[pg 48]{cS68} for definition), along with cohomological techniques to explore many nice basic properties of $\Aut_c(G)$ and $\Out_c(G)$ for a given finite group $G$.
We record a few of his results in Section 4.

After this, in 1980, H. Heineken \cite{hH80}, on the way to produce examples of finite groups in which all normal subgroups are characteristic, constructed finite $p$-groups $G$ with $\Aut(G) = \Aut_c(G)$. In particular, $\Out_c(G) \not= 1$ for all of these group $G$. His groups are $p$-groups of nilpotency class $2$. These groups also satisfy the property that $x\gamma_2(G) = x^G$ for all $x \in G - \gamma_2(G)$. 

Continuing in this direction, I. Malinowska \cite{iM92}, in 1992, constructed finite $p$-groups $G$ of nilpotency class $3$ and order $p^6$ such that $\Aut(G) = \Aut_c(G)$. In the same paper, she also constructed $p$-groups $G$ of nilpotency class $r$, for any prime $p > 5$ and any integer $r > 2$ such that $\Out_c(G) \not= 1$. 
In 1988, W. Feit and G. M. Seitz \cite[Section C]{FS88} proved that $\Out_c(G) = 1$ for all finite simple groups $G$.

Motivation for studying class preserving automorphisms of finite groups also arises from other branches of mathematics. 
T. Ono and H. Wada  proved that $\Out_c(G) = 1$ for the cases when $G$ is a free group, $SL_n(D)$, $GL_n(D)$ (where $D$ is an Euclidean domain), $S_n$ and  $A_n$ (the symmetric and alternating groups on $n$ symbols) in \cite{OW99, OW99a, hW99, hW00}. Their motivation for studying these things arose from ``Hasse principle" for smooth curves on number field. They associated a Shafarevich-Tate set to the given curve. The curve is said to enjoy ``Hasse principle" if the corresponding  
Shafarevich-Tate set is trivial. If the group $G$, involved in defining Shafarevich-Tate set, is finite then this set enjoys the group structure, which is isomorphic to $\Out_c(G)$. Interested reader can find all the details in \cite{tO98, tO01}.

The other motivation comes from integral groups rings.  Interested reader can refer \cite{mH04}. A. Hermann and Y. Li \cite{HL06} and M. Hertweck and E. Jespers \cite{HJ09} proved that $\Out_c(G) = 1$ for all Blackburn group. These groups were classified by N. Blackburn and satisfy the property that the intersection of all its non-normal subgroups is non-trivial. M. Hertweck \cite{mH01} constructed a class of Frobenius group $G$ such that $\Out_c(G) \not= 1$. He also proved that $\Out_c(G) = 1$ for all finite metabelian $A$-groups and for all $A$-groups having all Sylow subgroups elementary abelian.  By an $A$-group, we here mean a solvable group with abelian Sylow subgroups.

Other interesting examples of finite groups $G$ such that $\Out_c(G) \not= 1$ are given by F. Szechtman \cite{fS03}. He constructed a very concrete method to find such example. Interesting thing here is that his method produces examples in a variety of ways arising from different kind of Lie algebras, regular representations of finite fields and linear algebra. We would like to remark that his method produced examples of group having outer $n$-inner automorphims. An automorphism $\alpha$ of a group $G$ is said to be $n$-inner, if given any subset $S$ of $G$ with cardinality less than $n$, there exists an inner automorphism of $G$ agreeing with $\alpha$ on $S$. Notice that $2$-inner automorphisms are class preserving.

The present survey intersects with Section 5 of a survey article by I. Malinowska \cite{iM02}.

Our notation for objects associated with a finite multiplicative group $G$ is 
mostly standard.  The abelian group of all homomorphisms from an abelian group 
$H$ to an abelian group $K$ is denoted by $\Hom(H,K)$.
To say that some $H$ is a subset or a subgroup of $G$ we write 
$H \subseteq G$ or $H \leq G$ respectively. To indicate, in addition, that 
$H$ is properly contained in $G$, we write $H \subset G$, $H < G$ 
 respectively.  If $x,y \in G$, then $x^y$ denotes the conjugate element 
$y^{-1}xy \in G$ and $[x,y]$ denotes the commutator 
$x^{-1}y^{-1}xy = x^{-1}x^y \in G$. Here  $[x,G]$ 
denotes the set of all $[x,w]$, where $w$ runs over every element of $G$. Since $x^w = x[x,w]$, for all 
$w \in G$, we have $x^G = x[x,G]$. For $x \in G$, $\C_{H}(x)$ denotes the 
centralizer of $x$ in $H$, where $H \le G$. The center of $G$ will be 
denoted by $\Z(G)$. By $\gamma_2(G)$ we denote the commutator subgroup of a group $G$ which is generated by all the commutators $[x,y]$ in $G$.

\section{Some defintions and basic results}

We start with the concept of isoclinism of finite groups which was introduced by P. Hall \cite{pH40} (also see
\cite[pg. 39-40]{DY06} for details).

Let $X$ be a finite group and $\bar{X} = X/\Z(X)$. 
Then commutation in $X$ gives a well defined map
$a_{X} : \bar{X} \times \bar{X} \mapsto \gamma_{2}(X)$ such that
$a_{X}(x\Z(X), y\Z(X)) = [x,y]$ for $(x,y) \in X \times X$.
Two finite groups $G$ and $H$ are called \emph{isoclinic} if 
there exists an  isomorphism $\phi$ of the factor group
$\bar G = G/\Z(G)$ onto $\bar{H} = H/\Z(H)$, and an isomorphism $\theta$ of
the subgroup $\gamma_{2}(G)$ onto  $\gamma_{2}(H)$
such that the following diagram is commutative
\[
 \begin{CD}
   \bar G \times \bar G  @>a_G>> \gamma_{2}(G)\\
   @V{\phi\times\phi}VV        @VV{\theta}V\\
   \bar H \times \bar H @>a_H>> \gamma_{2}(H).
  \end{CD}
\]
The resulting pair $(\phi, \theta)$ is called an \emph{isoclinism} of $G$ 
onto $H$. Notice that isoclinism is an equivalence relation among finite 
groups. Each isoclinism class has a subgroup $G$ such that $\Z(G) \le \gamma_2(G)$, which is called a stem group of the family.

Let $G$ and $H$ be two isoclinic groups. Then the nilpotency class of $G$ is equal to the nilpotency class of $H$. Also the terms in the lower central series of $G$ are isomorphic to the respective terms of the lower central series of $H$. More precisely, one can say that the commutator structure of $G$ is similar to the commutator structure of $H$. The consequence of these simple results is the following theorem \cite[Theorem 4.1]{mY08}, which essentially says that the group $\Aut_c(G)$ is independent of the choice of the group $G$ in its isoclinism class.

\begin{thm}\label{thm1}
Let $G$ and $H$ be two finite non-abelian isoclinic groups. Then
$\Aut_c(G) \cong \Aut_c(H)$.
\end{thm}

Since every isoclinism class has a stem group, we readily get the following result.

\begin{cor}
It is sufficient to study $\Aut_c(G)$ for all finite groups $G$ such that $\Z(G) \le \gamma_2(G)$.
\end{cor}

Let $G$ be a finite group and  $N$ a non-trivial proper normal subgroup
of $G$. The pair $(G, N)$ is called a \emph{Camina pair} if $xN \subseteq x^G$ for all 
$x \in G-N$. A group $G$ is called a \emph{Camina group} if $(G, \gamma_2(G))$ is a Camina pair.
A basic result for Camina $p$-groups, which is due to R. Dark and C. M. Scoppola \cite{DS96} is the following.

\begin{thm}
Let $G$ be a finite Camina $p$-group. Then the nilpotency class of $G$ is at most $3$.
\end{thm}

A finite group $G$ is called \emph{purely non-abelian} if it does not have a non-trivial abelian direct factor.
An automorphism $\alpha$ of a group $G$ is called central if it commutes with every inner automorphism of $G$.
The group of all central automorphisms of a group $G$ is denoted by $\Autcent(G)$.
The following nice result is due to J. E. Adney and T. Yen \cite[Corollary 2]{AY65}.

\begin{thm}\label{aythm}
If $G$ is a purely non-abelian $p$-group, then $\Autcent(G)$ is also a $p$-group.
\end{thm}

\section{$\Out_c(G)$ for finite $p$-groups}

In this section we record several classes of finite $p$-groups $G$, for which it is decidable whether $\Out_c(G) = 1$ or not.  For the convenience of the reader we divide this section in several parts.

\vspace{.2in}

\noindent{\bf (i) $p$-groups of small orders.}  
\vspace{.2in}

In this part we consider $p$-groups of small orders and mention results on class-preserving automorphisms of $p$-groups of order at most $p^5$. Since $\Aut_c(G) = \Inn(G)$ for all abelian groups $G$, it is sufficient to consider the groups of order at least $p^3$.

For the proof of the following well known result, one can see \cite[pg 69]{mS86}.

\begin{prop}
Let $G$ be an extraspecial $p$-group.  Then $\Out_c(G) = 1$.
\end{prop}

Since every non-abelian group of order $p^3$ is extraspecial,  we get the following.

\begin{cor}
Let $G$ be a group of order $p^3$. Then $\Out_c(G) = 1$.
\end{cor}

We (jointly with L. R. Vermani) proved the following result in \cite{YV01}.

\begin{thm}
Let $G$ be a group of order $p^4$. Then $\Out_c(G) = 1$.
\end{thm}

Now we consider the groups of order $p^5$, where $p$ is an odd prime. 
Set 
\begin{equation}\label{eqn1a}
G_7 = \gen{a, b, x, y, z \;|\; {\mathcal R}, [b,a] = 1} \;\text{for}\; p \ge 3
\end{equation}
and 
\begin{equation}\label{eqn1b}
G_{10} = \gen{a, b, x, y, z \;|\; {\mathcal R}, [b,a] = x}\; \text{for}\;
p \ge 5,
\end{equation}

where ${\mathcal R}$ is defined as follows:
\begin{eqnarray*}
{\mathcal R} &=& \{a^p = b^p = x^p = y^p = z^p = 1\} \cup \{[x,y] = [x,z] = 
                 [y,z] = 1\}\\
          & & \quad \cup \{x^b = xz, y^b = y, z^b = z\} \cup \{x^a = xy, y^a =
                 yz, z^a = z\}.
\end{eqnarray*}

For $p = 3$ define a group $H$ of order $3^5$ by
\begin{eqnarray}\label{eqn1c}
H &=& \gen{a, b, c \;| \;a^3 = b^9 = c^9 = 1, [b,c] = c^3, [a,c] = b^3, [b, a] 
= c}\\
&=& \gen{a} \ltimes (\gen{b} \ltimes \gen{c}).\nonumber
\end{eqnarray}

\begin{remark}\label{rem2}
The group  $\phi_7(1^5)$ in the isoclinism family (7) of
\cite[Section 4.5]{rJ80} is isomorphic to $G_7$. The group $\phi_{10}(1^5)$ 
in the isoclinism family (10) of \cite[Section 4.5]{rJ80} is isomorphic to
$G_{10}$ for $p \ge 5$ and is isomorphic to $H$ for $p = 3$.
\end{remark}

The proof of the following theorem \cite[Theorem 5.5]{mY08} uses the classification of groups of order $p^5$ from \cite{rJ80} and Theorem \ref{thm1} above:
\begin{thm}\label{thm3a}
Let $G$ be a finite group of order $p^5$, where $p$ is an odd prime. Then
$\Out_c(G) \not= 1$ if and only if $G$ is isoclinic to one of the groups 
$G_7$, $G_{10}$ and $H$.
\end{thm}

This theorem, along with the examples of Wall \cite{gW47}, proves that $4$ is the smallest value of an integer $n$ such that all the groups of order less than or equal to $p^n$ has the property that $\Out_c(G) = 1$.

\vspace{.2in}

\noindent{\bf (ii)  Groups with large cyclic subgroups.}
\vspace{.2in}

In this part we study class preserving automorphisms of finite $p$-groups of order $p^n$ having cyclic subgroups of order $p^{n-1}$ or $p^{n-2}$.

The following result is Theorem 3.1 in \cite{YV00}
\begin{thm}
Let $G$ be a finite $p$-group having a maximal cyclic subgroup. Then $\Out_c(G) = 1$.
\end{thm}

For odd primes $p$, we have a similar result for finite $p$-groups $G$ having a cyclic subgroup of index $p^2$ \cite{YV02} and \cite{FN04}.
\begin{thm}
Let $G$ be a group of order $p^n$ having a cyclic subgroup of order $p^{n-2}$, where $p$ is an odd prime. Then $\Out_c(G) = 1$.
\end{thm}

But the situation is different for $p =2$. Let us set
\begin{eqnarray*}
G_1  &=& \langle x, y, z \;| \; x^{2^{m-2}} = 1 = y^{2} = z^{2},\; yxy = x^{1+2^{m-3}},
zyz = y,\\
& & zxz = x^{-1+ 2^{m-3}}\rangle;\\
G_2 &=& \langle x, y, z \;|\; x^{2^{m-2}} = 1 = y^{2} = z^{2},\;
yxy = x^{1+2^{m-3}},\; zxz = x^{-1+ 2^{m-3}},\\
& & zyz = yx^{2^{m-3}}\rangle.
\end{eqnarray*}

Then we get the following result \cite{YV02}:
\begin{thm}
Let $G$ be a group of order $2^n$, $n > 4$, such that $G$ has a cyclic normal subgroup of order $2^{n-2}$, but does not have any element of order $2^{n-1}$. Then $\Out_c(G) \not= 1$ if and only if $G$ is isomorphic to $G_1$ or $G_2$.
\end{thm}

\vspace{.2in}

\noindent{\bf (iii) Groups of nilpotency class $2$.}
\vspace{.2in}

Let $G$ be a finite nilpotent group of  class $2$. Let $\phi \in
\Aut_c(G)$. Then the map $g \mapsto g^{-1}\phi(g)$ is a homomorphism of
$G$ into $\gamma_2(G)$. This homomorphism sends $Z(G)$ to $1$. So it
induces a homomorphism $f_{\phi} \colon G/Z(G) \to \gamma_2(G)$, sending
$gZ(G)$ to $g^{-1}\phi(g)$, for any $g \in G$.  It is easily seen that
the map $\phi \mapsto f_{\phi}$ is a monomorphism of the group
$\Aut_c(G)$ into $\Hom(G/Z(G), \gamma_2(G))$.

Any $\phi \in \Aut_c(G)$ sends any $g \in G$ to some $\phi(g)
\in g^G$. Then $f_{\phi}(gZ(G)) = g^{-1}\phi(g)$ lies in $g^{-1}g^G =
[g,G]$.  Denote
\[  \{ \, f \in \Hom(G/Z(G), \gamma_2(G)) \mid f(gZ(G)) \in [g,G], \text{ for
all $g \in G$}\,\} \]
by $\Hom_c(G/Z(G), \gamma_2(G))$. Thus $f_{\phi} \in \Hom_c(G/Z(G),
\gamma_2(G))$ for all $\phi \in \Aut_c(G)$. On the other hand, if $f \in
\Hom_c(G/Z(G), \gamma_2(G))$, then the map sending any $g \in G$ to
$gf(gZ(G))$ is an automorphism $\phi \in \Aut_c(G)$ such that $f_{\phi}
= f$. Thus we have

\begin{prop}\label{prop1}
 Let $G$ be a finite nilpotent group of class 2. Then the
above map $\phi \mapsto f_{\phi}$ is an isomorphism of the group
$\Aut_c(G)$ onto $\Hom_c(G/$ $Z(G), \gamma_2(G))$.
\end{prop}

This correspondence gives the following result \cite[Theorem 3.5]{mY08}:

\begin{thm}\label{thm2}
Let $G$ be a finite $p$-group of class $2$. Let $\{x_1, x_2, \dots, x_d\}$ be
a minimal generating set for $G$ such that $[x_i, G]$ is cyclic, $1 \le i \le
d$. Then $\Out_c(G) = 1$.
\end{thm}

In particular, we get

\begin{cor}\label{cor1}
Let $G$ be a finite $p$-group of class $2$ such that $\gamma_2(G)$ is cyclic. 
Then $\Out_c(G) = 1$.
\end{cor}

\vspace{.2in}

\noindent{\bf (iv) Camina $p$-groups.}
\vspace{.2in}

Notice that extraspecial $p$-groups are Camina groups and $\Out_c(G) = 1$ for every extraspecial $p$-group $G$. The following result, which is proved in  \cite[Theorem 5.4]{mY07}, shows that these are the only Camina $p$-groups of class $2$ for which $\Out_c(G) = 1$.

\begin{thm}
Let $G$ be a finite Camina $p$-group of nilpotency class $2$. Then $|\Aut_c(G)| = |\gamma_2(G)|^d$, where $d$ is the number of elements in a minimal generating set for $G$.
\end{thm}

The following result \cite[Corollary 4.4]{mY07a} deals with finite $p$-groups $G$ such that $(G, \Z(G))$ is a Camina pair.

\begin{thm}
Let $G$ be a finite $p$-group of nilpotency class at least $3$ such that $(G,\Z(G))$ is a Camina pair and
$|\Z(G)| \ge p^2$. Then $|\Aut_c(G)| \ge |G|$.
\end{thm}

Since $(G, \Z(G))$ is a Camina pair for every Camina $p$-group $G$, we have the following immediate corollary:
\begin{cor}
Let $G$ be a finite Camina $p$-group of nilpotency class $3$ and $|\Z(G)| \ge p^2$. Then $|\Aut_c(G)| \ge |G|$.
\end{cor}

It is not diffucult to prove that $\Out_c(G) \not= 1$ for every Camina $p$-group of nilpotency class  $3$.

\section{Results of Sah}

In this section, we record a few results of C. H. Sah \cite{cS68}. This paper contains lots of reduction techniques and nice results, but it is very difficult to mention all the results here. For the readers, who are interested in class preserving automorphisms of group, we strongly recommend this beautiful article. 
 
\begin{thm}[Proposition 1.8, \cite{cS68}]
Let $G$ be a nilpotent group of class $c$. Then $\Aut_c(G)$ is a nilpotent group of class $c-1$.
\end{thm}

As an immediate consequence of this result one readily gets the following corollary.
\begin{cor}
Let $G$ be a nilpotent group of class $c$. Then $\Out_c(G)$ is a nilpotent group of class at most $c-1$.
\end{cor}

\begin{thm}[Theorem 2.9, \cite{cS68}]
Let $G$ be a finite solvable group. Then $\Aut_c(G)$ is solvable.
\end{thm}

\begin{thm}[Theorem 2.10, \cite{cS68}]
Let $G$ be a group admitting a composition series. Suppose that for each composition factor $F$ of $G$, the group $\Aut(F)/ \Inn(F)$ is solvable. Then $\Out_c(G)$ is solvable.
\end{thm}

Since every finite group admits a composition series and the Schreier's conjecture, i.e., the group of outer automorphisms of a finite simple group is solvable, holds true for all finite simple groups, the following result is a consequence of the above theorem.

\begin{cor}
Let $G$ be a  finite group. Then $\Out_c(G)$ is solvable.
\end{cor}

Sah also constructed examples of finite groups $G$ of order $p^{5n}$, $n \ge 3$, such that $\Out_c(G)$ is non-abelian. Thus, these examples contradict the intuitions of W. Burnside \cite[Note B]{wB11} that $\Out_c(G)$ should be abelian for all finite groups $G$.

\section{An upper bound for $\Aut_c(G)$ for a finite $p$-group $G$}

In this section, we list results from \cite{mY07}. Let $G$ be a finite $p$-group of order $p^n$.
Let $\{x_{1}, \cdots, x_{d}\}$ be any minimal generating set for $G$. 
Let $\alpha \in \Aut_{c}(G)$. Since $\alpha(x_{i}) \in x_{i}^G$ for 
$1 \leq i \leq d$, there are at the most $|x_{i}^{G}|$ choices for the image
of $x_i$ under $\alpha$. Thus it follows that
\begin{equation}
\label{bineq} |\Aut_{c}(G)| \leq \prod_{i=1}^{d} |x_{i}^G|.
\end{equation}

Let  $|\gamma_{2}(G)| = p^m$. Let $\Phi(G)$ denotes the Frattini subgroup 
of $G$. Since $\gamma_{2}(G)$ is contained in $\Phi(G)$, by the  Burnside 
basis theorem we have $d \leq n-m$. 
Notice that $|x_{i}^G| \leq |\gamma_{2}(G)| = p^m$ for all $i = 1, 2, 
\cdots, d$. So from \eqref{bineq} we get
\begin{equation}
\label{lbineq} |\Aut_{c}(G)| \leq p^{md} \leq (p^{m})^{n-m} = p^{m(n-m)}.
\end{equation}

\begin{thm}[Theorem 5.1, \cite{mY07}] \label{thm3}
Let $G$ be a finite $p$-group. Equality holds in \eqref{lbineq} if and only
if $G$ is either an abelian $p$-group, or a non-abelian  Camina special 
$p$-group.
\end{thm}

\begin{thm} [Theorem 5.5, \cite{mY07}] \label{thm4}
Let $G$ be a non-trivial $p$-group having order $p^{n}$. Then
\begin{equation}\label{gbineq}
|\Aut_{c}(G)| \leq
   \begin{cases}
     p^{\frac{(n^{2}-4)}{4}},  &\text{if $n$ is even;}\\
     p^{\frac{(n^{2}-1)}{4}}, &\text{if $n$ is odd.}
   \end{cases}
\end{equation}
\end{thm}

Since the bound in Theorem \ref{thm4} is attained by all non-abelian groups of order $p^3$  ($n$ odd) and the group constructed by W. Burnside ($n$ even), it follows that the bound in the theorem is optimal.

In the next theorem we are going to classify all finite $p$-groups which attain the bound in Theorem \ref{thm4}.
For the statement of the next result
we need the following group of order $p^6$, which is the group 
$\phi_{21}(1^6)$ in the isoclinism family ($21$) of \cite{rJ80}:
\begin{eqnarray} 
R &=& \langle \alpha, \alpha_1, \alpha_2, \beta, \beta_1, \beta_2| 
[\alpha_1, \alpha_2] = \beta, [\beta, \alpha_i]=\beta_i, 
[\alpha, \alpha_1]= \beta_2,\label{r}\\
& &\;\; [\alpha, \alpha_2] = \beta_1^{\nu},
\alpha^p=\beta^p=\beta_i^p = 1, \alpha_1^p=\beta_1^{(^p_3)},
\alpha_2^p=\beta_1^{-(^p_3)},\nonumber\\
& & \;\;i=1,2 \rangle,\nonumber
\end{eqnarray}
where $\nu$ is the smallest positive integer which is a non-quadratic residue
{$\mod p$} and $\beta_1$ and $\beta_2$ are central elements. We would like to remark here that this is not a minimal presentation of the group $G$. But it is sufficient for our purpose.

The following theorem gives a  classification of all finite $p$-groups $G$ for which the upper bound in Theorem \ref{thm4} is achieved by $|\Aut_c(G)|$. Recall that $W$ is the group of order $p^6$, defined in the second paragraph of the introduction.

\begin{thm}[Theorem 5.12, \cite{mY07}]\label{thm6}
Let $G$ be a non-abelian finite $p$-group of order $p^n$. Then
equality holds in \eqref{gbineq} if and only if one of the following holds: 
\begin{subequations}
\begin{align}
& \text{$G$ is an extra-special $p$-group of order $p^3$;}\\
& \text{$G$ is  a group of nilpotency class $3$ and order $p^4$;}\\
& \text{$G$ is a Camina special $p$-group isoclinic to the group $W$ and
    $|G|=p^6$;}\\
& \text{$G$ is isoclinic to $R$ and $|G|=p^6$.}
\end{align}
\end{subequations}
\end{thm}

\section{Some problems}

In this section, we formulate some research problems for finite $p$-groups, which certainly make sense for arbitrary finite groups. 

We know that $\Out_c(G) = 1$ for all finite simple group \cite{FS88}. But the proof uses classification of finite simple groups. We think that there should be a direct proof.

\begin{problem}
Let $G$ be a finite simple group. Without using classification of finite simple groups, prove that $\Out_c(G) = 1$.
\end{problem}

\begin{problem}
Let $G$ be a finite $p$-group of nilpotency class $2$. Find necessary and sufficient conditions on $G$ such that $\Out_c(G) = 1$.
\end{problem}

Let $G$ be a finite $p$-group of nilpotency class $2$. Then $\Aut_c(G) \le \Autcent(G)$.
\begin{problem}
Classify all finite $p$-groups $G$ of class $2$ such that $\Autcent(G) = \Aut_c(G)$.
\end{problem}

As we mentioned in the introduction that H. Heineken \cite{hH80} constructed examples of finite $p$-groups $G$ of nilpotency class $2$ such that $\Aut(G) = \Aut_c(G)$. This gives rise to the following natural problem.
\begin{problem}
Classify all finite $p$-groups $G$ of class $2$ such that $\Aut(G) = \Aut_c(G)$.
\end{problem}

The following problem arises from the work of Sah \cite[pg 61]{cS68}, which is also given in Malinowska's survey \cite{iM02} as Question 12.
\begin{problem}
Let $G$ be a finite $p$-group of nilpotency class $c > 2$. Give a sharp upper bound for the nilpotency class of $\Out_c(G)$.
\end{problem}

Let $G$ be a finite $p$-group with a minimal generating set $\{x_1, x_2, \cdots, x_d\}$. Then $|\Aut_c(G)| \le \Pi_{i=1}^d |x_i^G|$.

\begin{problem}
Classify all finite $p$-group $G$ with a minimal generating set $\{x_1, x_2,$ $\cdots, x_d\}$ such that 
 $|\Aut_c(G)| = \Pi_{i=1}^d |x_i^G|$.
\end{problem}

Since $|x^G| \le |\gamma_2(G)|$ for all $x \in G$, we can even formulate the following particular case of the preceeding problem.  
\begin{problem}\label{prob7}
Classify all finite $p$-groups $G$  such that $|\Aut_c(G)| = |\gamma_2(G)|^d$, where $d$ is the number of elements in a minimal generating set for $G$.
\end{problem}

Notice that  Problem \ref{prob7} has a solution (Theorem \ref{thm3}) when $\gamma_2(G) = \Phi(G)$.

Let $G$ be a finite $p$-group such that $\Z(G) \le \gamma_2(G)$, then it follows that $G$ is purely non-abelian and therefore $\Autcent(G)$ is a $p$-group (Theorem \ref{aythm}). 
\begin{problem}
Let $G$ be a finite $p$-group such that $\Out_c(G) \not= 1$ and $\Z(G) \le \gamma_2(G)$. Find a sharp lower bound for $|\Aut_c(G) \Autcent(G)|$ in term of $|G|$. More precisely, is it true that $|G| \le |\Aut_c(G) \Autcent(G)|$?
\end{problem}


\begin{thebibliography}{99}

\bibitem{AY65} 
J.~E.~Adney and T.~Yen, 
\emph{Automorphisms of a $p$-group}, 
Illinois J. Math. \textbf{9} (1965), 137-143.

\bibitem{wB11} W.~Burnside, 
\emph{Theory of groups of finite order}, 
2nd Ed. Dover Publications, Inc., 1955. Reprint of the 2nd edition 
(Cambridge, 1911).

\bibitem{wB13} W.~Burnside, 
\emph{On the outer automorphisms of a group}, 
Proc. London Math. Soc. (2)  \textbf{11} (1913), 40-42.

\bibitem{DY06}
E.~C.~Dade and M.~K.~Yadav,
\emph{Finite groups with many product conjugacy classes},
Israel Journal of Mathematics  \textbf{154} (2006), 29-49.

\bibitem{DS96}
R.~Dark and C.~M.~Scoppola,
\emph{On Camina groups of prime power order},
Journal of Algebra  \textbf{181} (1996), 787-802.

\bibitem{FS88}
W.~Feit and G.~M.~Seitz,
\emph{On finite rational groups and related topics},
Illinois J. Math., \textbf{33} (1988), 103-131.

\bibitem{FN04}
M.~Fuma and Y.~Ninomiya,
\emph{``Hasse principle'' for finite $p$-groups with cyclic subgroups of index
$p^2$}, Math. J. Okayama Univ.  \textbf{46} (2004), 31-38.

\bibitem{pH40}
P.~Hall,
\emph{The classification of prime power groups},
 Journal f\"ur die reine und angewandte Mathematik  \textbf{182} (1940), 
130-141.


\bibitem{hH80}
H.~Heineken,
\emph{Nilpotente Gruppen, deren s$\ddot{\text{a}}$mtliche Normalteiler 
charakteristisch sind}, Arch. Math. (Basel) \textbf{33} (1980), No. 6, 497-503.

\bibitem{HL06}
A.~Herman and Y.~Li,
\emph{Class preserving automorphisms of Blackburn groups},
J. Aust. Math. Soc. \textbf{80} (2006), 351-358.

\bibitem{mH01}
M.~Hertweck,
\emph{Class-preserving automorphisms of finite groups},
J. Algebra \textbf{241} (2001), 1-26.

\bibitem{mH04} 
M.~Hertweck,
\emph{Contributions to the integral representation theory of groups},
 Habilitationsschrift, University of Stuttgart (2004). Available at\\
 ``http://elib.uni-stuttgart.de/opus/volltexte/2004/1638''

\bibitem{HJ09}
M.~Hertweck and E.~Jespers,
\emph{Class-preserving automorphisms and normalizer property for Blackburn groups},
J. Group Theory \textbf{12} (2009), 157-169.



\bibitem{wG66}
W.~Gasch$\ddot{u}$tz,
\emph{Nichtabelsche p-gruppen besitzen $\ddot{a}$ussere p-automorphismen},
J.Algebra \textbf{4} (1966), 1-2.


\bibitem{rJ80}
R.~James,
\emph{The groups of order $p^6$ ($p$ an odd prime)},
Math. Comp. \textbf{34} (1980), 613-637.


\bibitem{YV00}
M.~Kumar and L.~R.~Vermani,
\emph{``Hasse principle'' for extraspecial $p$-groups},
Proc. Japan Acad.  \textbf{76}, Ser. A, No. 8 (2000), 123-125.

\bibitem{YV01}
M.~Kumar and L.~R.~Vermani,
\emph{``Hasse principle'' for groups of order $p^4$},
Proc. Japan Acad.  \textbf{77}, Ser. A, No. 6 (2001), 95-98.


\bibitem{YV02}
M.~Kumar and L.~R.~Vermani,
\emph{On automorphisms of some $p$-groups},
Proc. Japan Acad.  \textbf{78}, Ser. A, No. 4 (2002), 46-50.


\bibitem{iM92}
I.~Malinowska,
\emph{On quasi-inner automorphisms of a finite p-group}, 
Publ. Math. Debrecen \textbf{41} (1992), No. 1-2, 73-77.

\bibitem{iM02}
I.~Malinowska,
\emph{$p$-automorphisms of finite $p$-groups: problems and questions},
Advances in Group Theory (1992), 111-127.

\bibitem{tO98}
 T.~Ono,
\emph{A note on Shafarevich-Tate sets for finite groups},
Proc. Japan Acad. \textbf{74} Ser. A (1998), 77-79.

\bibitem{tO01}  
T.~Ono,
\emph{On Shafarevich-Tate sets},
Advanced Studies in Pure Mathematics: Class Field Theory - Its Centenary and prospect \textbf{30} (2001),
537-547.

\bibitem{OW99}
T.~Ono and H.~Wada,
\emph{``Hasse principle'' for free groups},
Proc. Japan Acad., \textbf{75} Ser.  A (1999), 1-2.

\bibitem{OW99a}
T.~Ono and H.~Wada,
\emph{``Hasse principle'' for symmetric and alternating groups},
 Proc. Japan Acad., \textbf{75} Ser. A (1999), 61-62.


\bibitem{cS68}
C.~H.~Sah,
\emph{Automorphisms of finite groups}, Journal of Algebra \textbf{10} (1968),
47-68.

\bibitem{mS86} 
M.~Suzuki,
\emph{Group Theory II},
Springer, New York - Berlin - Heidelberg - Tokyo (1986).

\bibitem{fS03}
F.~Szechtman,
\emph{$n$-inner automorphisms of finite groups}, Proc. Amer. Math. Soc.
 \textbf{131} (2003), 3657-3664.

\bibitem{hW99}
H.~Wada,
\emph{``Hasse principle'' for $SL_{n}(D)$},
Proc. Japan Acad. \textbf{75} Ser. A (1999), 67-69.

\bibitem{hW00}
H.~Wada,
\emph{``Hasse principle'' for $GL_{n}(D)$},
Proc. Japan Acad. \textbf{76} Ser. A (2000), 44-46.


\bibitem{gW47}
G.~E.~Wall,
\emph{Finite groups with class preserving outer automorphisms},
J. London Math. Soc.  \textbf{22} (1947), 315-320.

\bibitem{mY07}
M.~K.~Yadav,
\emph{Class preserving automorphisms of finite $p$-groups}, 
J. London Math. Soc.  \textbf{75(3)} (2007), 755-772.


\bibitem{mY07a}
M.~K.~Yadav,
\emph{On automorphisms of finite $p$-groups}, 
J. Group Theory \textbf{10} (2007), 859-866.

\bibitem{mY08}
M.~K.~Yadav,
\emph{On automorphisms of some finite $p$-groups}, 
Proc. Indian Acad. Sci. (Math. Sci.) \textbf{118(1)} (2008), 1-11.


\end{thebibliography}
\end{document}